\documentclass[12pt,leqno]{amsart}
\usepackage{amssymb,verbatim,enumerate,ifthen,amsmath}
\usepackage[mathscr]{eucal}
\usepackage[utf8]{inputenc}
\usepackage[T1]{fontenc}
\usepackage{cite}
\def\N{\mathbb{N}}
\def\R{\mathbb{R}}

\def\L{\mathbb{L}}
\def\Z{\mathbb{Z}}

\def\M{\mathscr{M}}

\def\H{\mathscr{H}}

\newtheorem{theorem}{Theorem}[section]
\newtheorem*{theorem*}{Theorem}
\def\Thm#1#2{\ifthenelse{\equal{#1}{*}}{\begin{theorem*}#2\end{theorem*}}
             {\begin{theorem}\label{T#1}#2\end{theorem}}}
\newtheorem{Atheorem}{Theorem}

\def\thm#1{Theorem~\ref{T#1}}
\newtheorem{proposition}[theorem]{Proposition}
\newtheorem*{proposition*}{Proposition}
\def\Prp#1#2{\ifthenelse{\equal{#1}{*}}{\begin{proposition*}#2\end{proposition*}}
{\begin{proposition}\label{P#1}#2\end{proposition}}}
\def\prp#1{Proposition~\ref{P#1}}

\newtheorem{corollary}[theorem]{Corollary}
\newtheorem*{corollary*}{Corollary}
\def\Cor#1#2{\ifthenelse{\equal{#1}{*}}{\begin{corollary*}#2\end{corollary*}}
             {\begin{corollary}\label{C#1}#2\end{corollary}}}

\newtheorem{lemma}[theorem]{Lemma}
\newtheorem*{lemma*}{Lemma}
\def\Lem#1#2{\ifthenelse{\equal{#1}{*}}{\begin{lemma*}#2\end{lemma*}}
             {\begin{lemma}\label{L#1}#2\end{lemma}}}
\def\lem#1{Lemma~\ref{L#1}}

\theoremstyle{definition}
\newtheorem{remark}[theorem]{Remark}
\newtheorem*{remark*}{Remark}
\def\Rem#1#2{\ifthenelse{\equal{#1}{*}}{\begin{remark}\rm #2\end{remark}}
             {\begin{remark}\label{R#1}\rm #2\end{remark}}}

\newtheorem{example}[theorem]{Example}
\newtheorem*{example*}{Example}
\def\Exa#1#2{\ifthenelse{\equal{#1}{*}}{\begin{example*}\rm #2\end{example*}}
             {\begin{example}\label{Ex#1}\rm #2\end{example}}}

\def\eq#1{{\rm(\ref{E#1})}}
\def\Eq#1#2{\ifthenelse{\equal{#1}{*}}
  {\begin{equation*}\begin{aligned}#2\end{aligned}\end{equation*}}
  {\begin{equation}\begin{aligned}\label{E#1}#2\end{aligned}\end{equation}}}

\begin{document}
\begin{flushright}
\end{flushright}
\vspace{5mm}

\date{\today}

\title{Revisiting 
Fekete's Lemma, Subadditive and Periodic Sequences\\
}

\author[A. R. Goswami]{Angshuman R. Goswami}
\address{Department of Mathematics, University of Pannonia, 
Veszprém, Egyetem u. 10, H-8200.}
\email{goswami.angshuman.robin@mik.uni-pannon.hu}

\subjclass[2000]{Primary 39B62; Secondary 26A48, 52A30,40-00}
\keywords{ Subadditivite sequence, Subadditivite functions and Fekete's lemma }

\begin{abstract}
In this paper, we present an alternative proof of Fekete's Lemma. We demonstrate that for any subadditive sequence, it is possible to construct a subadditive function that exactly interpolates the sequence. Using this result, along with Hille's theorem on subadditive functions, we naturally arrive at Fekete's Lemma. Additionally, we provide an explicit formula for determining the largest subadditive minorant of a given sequence. We explore a sandwich-type result and derive a discrete version of the Hyers-Ulam type stability theorem. For approximately periodic sequences, we offer a decomposition result. In the final section, we propose two characterization theorems for ordinary periodic sequences.\\

The motivation, research progress, and other important details are covered in the introduction.
\end{abstract}
\maketitle
\section*{Introduction}
Throughout this paper, the symbols $\N$, $\Z$, and $\R$ denote the sets of natural, integer, and real numbers respectively. Primarily, this paper can be subdivided into two folds. In the first section we discuss subadditive sequences, while in the second part, we investigate periodic sequences.\\

Subadditivity is a natural phenomenon. It can be seen in several branches of mathematics in multiple formats. Such as modulus, norm, probability addition, triangle side rules, etc. A sequence $\left<u_i\right>_{i=1}^{\infty}$ is said to be subadditive if the following discrete functional inequality is satisfied
\Eq{654}{
u_{m+n}\leq u_m+u_n \qquad\mbox{for all} \qquad m,n\in\N.
}
The class of all subadditive sequences is a convex cone as it is closed under addition and multiplication with non-negative scalars. Unlike monotone, convex, and arithmetic progressions; ordering plays the decisive role in determining the existence of subadditive property in a sequence. For example  
$\big<4,1,0,-1,-4,-2,-3\big>$ is a subadditive sequence when considering the initial term as $u_1=4$. But, the very same sequence can not be labeled as a subadditive if we start marking its elements from $0$. The same characteristic irregularity can be observed for the sequence $\big<1,-1,0\big>$ due to the separate selection of discrete domains. To maintain technical preciseness and to avoid complications; we can also state that a sequence $\left<u_i\right>_{i=0}^{\infty}$ with $u_0=0$ is also subadditive if the inequality mentioned in \eq{654} holds.\\

Many research articles have been published that deal with subadditive functions and sequential subadditivity. The study of approximately subadditive and higher order subadditive functions, discovering the linkage of subadditive function with convex and periodic functions, several important limit theorems have been investigated in the last few decades. The detailed research can be found in the papers \cite{Fekete,Hillee,Hille,Schechter,Matkowski, Rosenbaum,Matkowskii,Goswamii,Burai,Arpad} and the reference therein. 

The main objective of our paper is to provide an alternative proof of Fekete's lemma, which can be stated as follows:  
"If  $\left<u_n\right>_{n=1}^{\infty}$ be a 
subadditive sequence, then
$$\lim_{n\to\infty}\dfrac{u_n}{n}=\underset{n\geq1}{\inf}\dfrac{u_n}{n}\quad \quad \quad (n\in\N)."$$
Fekete proved this result in 1923 in his paper  \cite{Fekete}. Around a quarter-century later, Hille came up with the functional version of Fekete's subadditive lemma. The proofs of both of these theorem are analogous. We establish a bridge between these two well-known theorems showing that Hille's theorem together with some other results yields Fekete's lemma.\\

By connecting the dots of a monotone sequence graph, we naturally obtain a monotone function, and the same holds for convex and arithmetic sequences. However, this is not true for many other types of sequences. For example, in geometric and Fibonacci sequences, simply joining the discrete points with line segments does not produce a function that preserves all the properties of the original sequence. In most cases, such functions require precise, rigorous formulations. The Gamma function derived from the factorial sequence is a notable example of it. In this context, the class subadditive sequences share similarities with monotone and convex sequences, but validating this requires careful consideration of multiple cases.
In other words, for a sequence $\left<u_i\right>_{i=1}^{n}$ possessing subadditivity, we can define a continuous subadditive function $f:[1,\infty[\to\R$ that perfectly interpolates the sequence. \\

In mathematics, the terminology of periodicity was known a long time ago. The trigonometric functions, Dirichlet functions, and most of the wave functions are just some of the well-known examples of it. The theory of entire Fourier analysis, wavelet analysis, and many other branches of mathematics are built on periodic functions. Since sequential periodicity has immense applications in computer science and many other fields; many researchers are actively working in this direction. Some of the most important results can be found in \cite{Gilbert, Heuberger, Muchnik, Kurshan, Khamidullin} and the references therein. \\

The periodic sequences carry several nice properties. For instance if $\left<u_{_n}^{^{1}}\right>_{n=1}^{\infty},\cdots,\mbox{and} \left<u_{_n}^{^{k}}\right>_{n=1}^{\infty}$ are $k$ periodic sequences with length $\L$ and $a_1,\cdots a_k\in\R$, then the sequence $\left<\overset{k}{\underset{i=1}{\sum}}a_{_i}u_{_n}^{^{i}}\right>_{n=1}^{\infty}$ also possesses periodicity with same periodic length. This yields that the class of periodic sequences with length $\L$ is a centrally symmetric periodic cone. Moreover, this class can also be treated as a commutative group under the operation point-wise sequence addition. \\ 

Several approximate notions of periodicity have been proposed. Out of all these, the eventually periodic and asymptotically periodic sequences are broadly studied. We found that both of these periodic sequence classes can be studied under the umbrella of $\varepsilon$-periodicity.  
A sequence $\left<u_n\right>_{n=1}^{\infty}$ is said to be $\varepsilon$-periodic with the period $\L$ ($\L\in\N$) if for any $i\in\N$, the following two discrete functional inequalities are satisfied simultaneously
\Eq{*}{
u_i-u_{_{i+k\L}}\leq \varepsilon \quad\mbox{and}\quad u_{_{i+k\L}}-u_i\leq \varepsilon\quad\quad\mbox{for all}\quad k\in\Z \quad\mbox{such that}\quad i+k\L\in\N.
}
We establish that a $\varepsilon$-periodic sequence can be expressed as the algebraic summation of a periodic and a controlled sequence that fluctuates within the interval $\left[-\dfrac{\varepsilon}{2},\dfrac{\varepsilon}{2}\right].$ The reverse implication of this statement is also validated. In function theory, such decomposition results are termed as stability theorems. Some interesting details in this direction can be found in the papers \cite{Hyers,Ulam}.\\
  
Besides, we propose two characterizations of periodic sequences. We prove that a periodic sequence with periodicity $\L$ can be partitioned into a maximum $\L$ number of non-intersecting constant sub-sequences.  Also we show that for a periodic sequence $\left<u_n\right>_{n=1}^{\infty}$ with period $\L$, we can derive a continuous, periodic function $f:[1,\infty[\to\R,$ with same periodicity such that $f(n)=n$ holds for all $n\in\N$.
\section{On Subadditive sequences }
We start our investigation with some basic observations. The proofs of these statements are not difficult and hence only the statements are mentioned.
\Prp{0}{
Let $\left<u_n\right>_{n=1}^{\infty}$ be 
an arbitrary sequence then depending on the characterization; the following statements hold.
\begin{enumerate}[(i)]
\item If $\left<u_n\right>_{n=1}^{\infty}$ is non-negative and decreasing then the sequence also possesses subadditivity. 
\\
\item If $\left<u_n\right>_{n=1}^{\infty}$  possesses subadditivity satisfying $u_1\leq 0$; then the sequence is non-positive.\\
\item If $\left<u_n\right>_{n=1}^{\infty}$  is subadditive then for any $n>m$, the inequality $u_n-u_m\leq u_{n-m}$ always holds.
\end{enumerate}
}

The two propositions below together prove that for any given sequence $\left<u_n\right>_{n=1}^{\infty}$, there exists a unique 
subadditive majorant $\left<v_n\right>_{n=1}^{\infty}$ bounded above by $\left<u_n\right>_{n=1}^{\infty}$. In this first proposition, we will see the construction of such a subadditive sequence. 

\Prp{10}{Let $\left<u_n\right>_{n=1}^{\infty}$ be any arbitrary sequence. Then there exists a subadditive 
\newline sequence $\left<v_n\right>_{n=1}^{\infty}$ such that the expression $v_n\leq u_n$ holds for all $n\in\N$.}
\begin{proof}
To prove this proposition, we assume $n\in\N$ and $n_1,\cdots,n_k$ be any arbitrary partition of it. In other words $n_1,\cdots,n_k\in\N$ such that    $n_1+\cdots+n_k=n$. We define the sequence $\left<v_n\right>_{n=1}^{\infty}$ as follows
\Eq{22}{
v_n:=\min\Big\{u_{n_1}+\cdots+u_{n_k}\quad \big|\quad  n_1,\cdots,n_k\in\N  \quad \mbox{satisfying}\quad n_1+\cdots+n_k=n\Big\}
.}
By the construction of $\left<v_n\right>_{n=1}^{\infty}$, it is clear that  $v_n\leq u_n$ holds for all $n\in\N$. We only need to show that the sequence $\left<v_n\right>_{n=1}^{\infty}$ possesses subadditivity. We consider $m,n\in\N$. Then there exists respective partitions of $m$ and $n$ such that $m=m_1+\cdots m_l$ and $n=n_1+\cdots+n_k$ with $m_1,\cdots m_l,\,n_1,\cdots, n_k\in\N$ satisfying the following two inequalities
\Eq{33}{
v_m=u_{m_1}+\cdots+u_{m_l}\quad\mbox{and}\quad
v_n=u_{n_1}+\cdots+u_{n_k}.}
The partitions of $m$ and $n$ provide a partition for $m+n$ as well which can be represented as    
$
m+n= m_1+\cdots m_l+ n_1+\cdots+ n_k$ with $ m_1,\cdots m_l$ and $n_1,\cdots, n_k\in\N.
$
From the construction of the sequence  $\left<v_n\right>_{n=1}^{\infty}$ and \eq{33}, we can conclude the following inequality
\Eq{*}{
v_{m+n}\leq u_{m_1}+\cdots+u_{m_l}+ u_{n_1}+\cdots+u_{n_k}= v_m+v_n.
}
This yields that the sequence $\left<v_n\right>_{n=1}^{\infty}$ is subadditive and proves our assertion.\\
\end{proof}
In the next proposition, we will prove the uniqueness part, that we have discussed before.
\Prp{86}{
The sequence $\left<v_n\right>_{n=1}^{\infty}$ as defined in \eq{22} is the largest subadditive majorant bounded above by $\left<u_n\right>_{n=1}^{\infty}$.
} 
\begin{proof}
Subadditivity of $\left<v_n\right>_{n=1}^{\infty}$ was already established in the previous proposition.
\newline
We assume there  exists another subadditive  sequence $\left<w_n\right>_{n=1}^{\infty}$ satisfying the inequality
\Eq{77}{v_n\leq w_n\leq u_n \quad\mbox {for all} \quad n\in\N.}
Let $n\in\N$ be arbitrary such that $n_1,\cdots,n_k$ is a partition of it. Then the inequality below is satisfied
\Eq{*}{w_n\leq w_{n_1}+\cdots w_{n_k}\leq u_{n_1}+\cdots u_{n_k}.}
After taking the minimum of the rightmost expression of the above inequality concerning all possible partitions of $n$, we arrive at $w_n\leq v_n$. Since $n$ is arbitrary, we have
\Eq{*}{w_n\leq v_n\leq u_n \quad\mbox {for all} \quad n\in\N.}
This inequality together with \eq{77} yields the uniqueness of the subadditive sequence $\left<v_n\right>_{n=1}^{\infty}$ and completes the proof.
\end{proof}
Based on \prp{10} we can formulate the following sandwich-type result as an analogy.
\Cor{78}{Let $\left<w_n\right>_{n=1}^{\infty}$ and $\left<u_n\right>_{n=1}^{\infty}$ are two sequences satisfying the discrete functional inequality
$w_n\leq u_{n_1}+\cdots u_{n_k}$ for all $n\in\N$; where $n_1,\cdots, n_k\in\N$ is an arbitrary partition of $ n.$ Then there exists a subadditive sequence $\left<v_n\right>_{n=1}^{\infty}$ such that $w_n\leq v_n\leq u_n$ holds for all $n\in\N.$}
\begin{proof}
We can construct $\left<v_n\right>_{n=1}^{\infty}$ as in \eq{22} and from there we can validate the result.
\end{proof}
\section{On Fekete's Lemma}
A function $f:I(\subseteq\R)\to\R$ is said to be subadditive if for any $x,y\in I$ with $x+y\in I$ the following functional inequality is satisfied
\Eq{10,000}{
f(x+y)\leq f(x)+f(y).
}
The proposition below shows that for any sequence possessing subadditivity, there must be an underlying continuous subadditive function that perfectly interpolates the sequence. 
\Thm{2200}{Let $\big<u_n\big>_{n=1}^{\infty}$
be a subadditive sequence. Then there must exist a continuous, almost everywhere 
differentiable subadditive function $f:[1,\infty[\to\R$ such that $f(n)=u_n$ holds for all $n\in\N.$ On the other hand, if $f:[1,\infty[\to\R$ is a subadditive function then $\big<f(n)\big>_{n=1}^{\infty}$ represents a subadditive sequence.
}
\begin{proof}
To establish the statement, we define the function $f:[1,\infty[\to\R$ as follows
\Eq{4000}{
f(u):=tu_n+(1-t)u_{n+1}\quad\mbox{where}\quad u:=tn+(1-t)(n+1)\quad \Big(t\in[0,1]\quad\mbox{and}\quad n\in\N\Big).
}
From the construction of $f$, continuity is obvious. Moreover, $f$ is differentiable everywhere in $[1,\infty[$ except the set $\N$. Clearly, $f(n)=u_n$ holds for all $n\in\N$.\\

Therefore, we are only left to show the subadditivity of the function $f$.
We consider two arbitrary elements $x\in[n_1,n_{1}+1]$ and $y\in[n_2,n_{2}+1]$ such that
\Eq{*}{
x:=t_1{n_{1}}+(1-t_1)({n_{1}}+1)\,\,\mbox{and}\,\, y:=t_2{n_{2}}+(1-t_2)({n_{2}+1})\,\, \mbox{where}\,\, n_1,n_2 \in\N, \,\, t_1,t_2\in [0,1].
}
This yields two possibilities. At the beginning, we will consider $x+y\in[n_1+n_2,n_1+n_2+1]$. While in the second case, we will go through the scenario when $x+y\in[n_1+n_2+1,n_1+n_2+2].$\\

At first, we assume $x+y:=t(n_1+n_2)+(1-t)(n_1+n_2+1)$; where $t\in[0,1]$. This provides $t=t_1+t_2-1.$ To show the subadditivity of the function $f:[1,\infty[\to\R$, by utilizing the sequential subadditivity of $\big<u_n\big>_{n=1}^{\infty}$, we can proceed as follows
\Eq{5000}{
f(x)+f(y)-f(x+y)&=t_1u_{_{n_1}}+(1-t_1)u_{_{n_1+1}}+t_2u_{_{n_2}}+(1-t_2)u_{_{n_2+1}}\\
&\qquad\qquad\qquad\qquad\qquad\qquad\qquad\qquad\quad-tu_{_{n_1+{n_2}}}-(1-t)u_{_{n_1+{n_2}+1}}\\
&=t_1\Big(u_{_{n_1}}+u_{_{n_1+n_2+1}}-u_{_{n_1+1}}-u_{_{n_1+n_2}}\Big)+t_2\Big(u_{_{n_2}}+u_{_{n_1+n_2+1}}-\\
&\qquad \qquad\quad\,\,\,  u_{_{n_2+1}}-u_{_{n_1+n_2}}\Big)+\Big(u_{_{n_1+1}}+u_{_{n_2+1}}+u_{_{n_1+n_2}}-2u_{_{n_1+n_2+1}}\Big)\\
&= t_1\Big(u_{_{n_1}}+u_{_{n_1+n_2+1}}-u_{_{n_1+1}}-u_{_{n_1+n_2}}\Big)+t_2\Big(u_{_{n_2}}+u_{_{n_1+n_2+1}}\\
&\quad \quad  u_{_{n_2+1}}-u_{_{n_1+n_2}}\Big)+(t_1+t_2)\Big(u_{_{n_1+1}}+u_{_{n_2+1}}+u_{_{n_1+n_2}}-2u_{_{n_1+n_2+1}}\Big)\\
&\qquad\qquad\qquad \,\,\,\,\,+(1-t_1-t_2)\Big(u_{_{n_1+1}}+u_{_{n_2+1}}+u_{_{n_1+n_2}}-2u_{_{n_1+n_2+1}}\Big)\\
&=t_1\Big(u_{_{n_1}}+u_{_{n_2+1}}-u_{_{n_1+n_2+1}}\Big)+t_2\Big(u_{_{n_2}}+u_{_{n_1+1}}-u_{_{n_1+n_2+1}}\Big)\\
&\qquad\qquad\qquad\qquad (1-t_1-t_2)\Big(u_{_{n_1+1}}+u_{_{n_2+1}}+u_{_{n_1+n_2}}-2u_{_{n_1+n_2+1}}\Big)\\
&\geq \min\{t_1,t_2\}\Big(\big(u_{_{n_1}}+u_{_{n_2+1}}-u_{_{n_1+n_2+1}}\big)+\big(u_{_{n_2}}+u_{_{n_1+1}}-u_{_{n_1+n_2+1}}\big)\Big)\\
&\qquad\qquad\qquad\,\,\,\,\, +(1-t_1-t_2)\Big(u_{_{n_1+1}}+u_{_{n_2+1}}+u_{_{n_1+n_2}}-2u_{_{n_1+n_2+1}}\Big).
} 
In the rightmost side of the above expression, the terms, $u_{_{n_1}}+u_{_{n_2+1}}-u_{_{n_1+n_2+1}}$ and 
\newline 
$u_{_{n_2}}+u_{_{n_1+1}}-u_{_{n_1+n_2+1}}$ are non-negative. While $1-t_1-t_2$ (or $-t$) is non-positive.\\
 
If $u_{_{n_1+1}}+u_{_{n_2+1}}+u_{_{n_1+n_2}}-2u_{_{n_1+n_2+1}}\leq 0$; then the above inequality  implies
\eq{10,000}.\\

On the other hand, if $u_{_{n_1+1}}+u_{_{n_2+1}}+u_{_{n_1+n_2}}-2u_{_{n_1+n_2+1}}\geq 0$; we can extend \eq{5000} as follows
\Eq{*}{
f(x)+f(y)-f(x+y)&\geq\min\{t_1,t_2\}\Big(u_{_{n_1}}+u_{_{n_1+1}}+u_{_{n_2}}+u_{_{n_2+1}}-2u_{_{n_1+n_2+1}}\Big)\\
&\quad\qquad\qquad\qquad +(1-t_1-t_2)\Big(u_{_{n_1+1}}+u_{_{n_2+1}}+u_{_{n_1+n_2}}-2u_{_{n_1+n_2+1}}\Big)\\
&\geq\big(\min\{t_1,t_2\}+(1-t_1-t_2)\big)\Big(u_{_{n_1+1}}+u_{_{n_2+1}}+u_{_{n_1+n_2}}-2u_{_{n_1+n_2+1}}\Big)\\
&=\big(1-\max\{t_1,t_2\}\big)\Big(u_{_{n_1+1}}+u_{_{n_2+1}}+u_{_{n_1+n_2}}-2u_{_{n_1+n_2+1}}\Big)\\
&\geq 0.}
Same conclusions for the two sub-cases yields subadditivity of $f$; for $x+y\in[n_1+n_2,n_1+n_2+1]$.\\

Now we let, $x+y:=t(n_1+n_2+1)+(1-t)(n_1+n_2+2)$; where $t\in[0,1]$. This gives $t=t_1+t_2.$ To show the subadditivity of the function $f:[1,\infty[\to\R$, we proceed as follows
\Eq{*}{
f(x)+f(y)-f(x+y)&=t_1u_{_{n_1}}+(1-t_1)u_{_{n_1+1}}+t_2u_{_{n_2}}+(1-t_2)u_{_{n_2+1}}\\
&\qquad\qquad\qquad\qquad\qquad\qquad\qquad\,
-tu_{_{n_1+n_2+1}}-(1-t)u_{_{n_1+n_2+2}}\\
&=t_1\big(u_{_{n_1}}-u_{_{n_1+1}}\big)+t_2\big(u_{_{n_2}}-u_{_{n_2+1}}\big)+u_{_{n_1+1}}+u_{_{n_2+1}}-u_{_{n_1+n_2+2}}\\
&\qquad\qquad\qquad\qquad\qquad\qquad\quad\,\,\,\,
+(t_1+t_2)\big(u_{_{n_1+n_2+2}}-u_{_{n_1+n_2+1}}\big)\\
&\geq t_1\big(u_{_{n_1}}+u_{_{n_1+n_2+2}}-u_{_{n_1+1}}-u_{_{n_1+n_2+1}}\big)+t_2\big(u_{_{n_2}}+u_{_{n_1+n_2+2}}\\
&\qquad  \quad -u_{_{n_2+1}}-u_{_{n_1+n_2+1}}\big)+(t_1+t_2)\big(u_{_{n_1+1}}+u_{_{n_2+1}}-u_{_{n_1+n_2+2}}\big)\\
&=t_1\Big(u_{_{n_1}}+u_{_{n_2+1}}-u_{_{n_1+n_2+1}}\Big)+t_2\big(u_{_{n_1}}+u_{_{n_1+1}}-u_{_{n_2}}-u_{_{n_1+n_2+1}}\big)
.}
Since all the terms in the rightmost part of the above inequality are non-negative, we get 
\eq{10,000} whenever $x+y\in[n_1+n_2+1,n_1+n_2+2]$ as well. Both the cases together prove the statement and it validates the first part of the theorem.\\

To show the converse part, let $f:[1,\infty[\to\R$ be a subadditive function. Then for any $x,y\in[1,\infty[$, it will satisfy the inequality $\eq{10,000}$. We assume $n\in\N$ to be arbitrary. Then for any $n_1,n_2\in\N$ with $n_1+n_2\in\N$, we can conclude the following inequality
$$f(n)=f(n_1+n_2)\leq f(n_1+f(n_2).$$
This shows that the $\left<f(n)\right>_{n=1}^{\infty}$ possesses sequential subadditivity and completes our proof. 
\end{proof}
The above theorem also shows that subadditive sequences are just the discrete versions of subadditive functions. To establish Fekete's lemma we will need the following result. 
\Lem{991}{Let $\left<u_n\right>_{n=1}^{\infty}$ be any sequence. Then for the function $f:[1,\infty[\to\R$ defined in \eq{4000}, the following equality holds
\Eq{771}{
\inf_{u\in[1,\infty[}\dfrac{f(u)}{u}=\inf_{n\in\N}\dfrac{u_n}{n}.
}
}
\begin{proof} 
The inequality $\underset{u\in[1,\infty[}{\inf}\dfrac{f(u)}{u}\leq\underset{n\in\N}{\inf}\dfrac{u_n}{n}$ is obvious. Therefore, to prove this lemma, it will be sufficient to show the reversed inequality. For that, we need the following result\\

"For $a_1,a_2\in\R$ and $b_1,b_2>0$, the following inequality is satisfied
\Eq{22111}{
\min\bigg\{\dfrac{a_1}{b_1},\dfrac{a_2}{b_2}\bigg\}\leq \dfrac{a_1+a_2}{b_1+b_2}\leq \max\bigg\{\dfrac{a_1}{b_1},\dfrac{a_2}{b_2}\bigg\}^".
}
The proof of the above statement and a generalized form of it can be found in one of our recently submitted papers \cite{Goswami}.
Now for any $u\in[1,\infty[$ with $u:=tn+(1-t)(n+1)$ ($n\in\N$); we have $\dfrac{f(u)}{u}=\dfrac{tu_n+(1-t)u_{n+1}}{tn+(1-t)(n+1)}$.
 Then by using \eq{22111}, we can compute the following inequality
\Eq{*}{
\min\bigg\{\dfrac{u_n}{n},\dfrac{u_n+1}{n+1}\bigg\}\leq\dfrac{f(u)}{u}\leq\max\bigg\{\dfrac{u_n}{n},\dfrac{u_n+1}{n+1}\bigg\}.
}
Since $u$ is arbitrary, by using the leftmost inequality we can conclude, 
$\underset{u\in[1,\infty[}{\inf}\dfrac{f(u)}{u}\geq\underset{n\in\N}{\inf}\dfrac{u_n}{n}.$ 
\newline And this completes the proof.
\end{proof}
The following result was proposed by 
Swedish-American Mathematician Einar Hille and it is known as Hille's subadditive theorem. 
\Thm{999}{\normalfont{\textbf{ [Hille's subadditive theorem] }}   For every subadditive function $f:[1,\infty[\to\R$, the limit $\underset{t\to\infty}{\lim}\dfrac{f(u)}{u}$ exits and is equal to $\underset{u\in[1,\infty[}{\inf}\dfrac{f(u)}{u}.$}
The upcoming result is known as Fekete's subadditive lemma. The classical proof of it includes concepts such as $limsup$, $liminf$, and tools like division algorithm. Through our proof, we establish the proper linkage of  Hille's theorem and Fekete's lemma on subadditivity. 
\Cor{1000} {\normalfont{\textbf{ [Fekete's subadditive lemma] }}For every  subadditive sequence $\left<u_n\right>_{n=1}^{\infty}$, the limit $\underset{n\to\infty}{\lim}\dfrac{u_n}{n}$ exits and is equal to infimum $\underset{n\in\N}{\inf}\dfrac{u_n}{n}.$ }
\begin{proof}
By using the interpolations of the function $f$ in \thm{2200}, Hille's subadditive theorem (\thm{999}), and equality \eq{771} of \lem{991} we can compute the following 
\Eq{*}{
\underset{n\to\infty}{\lim}\dfrac{u_n}{n}=\lim_{n\to\infty}\dfrac{f(n)}{n}=\lim_{u\to\infty}\dfrac{f(u)}{u}=\inf_{u\in[1,\infty[}\dfrac{f(u)}{u}=\inf_{n\in\N}\dfrac{u_n}{n}.
} 
This completes the proof.\\
\end{proof}
This establishment also shows that Fekete's subadditive lemma is actually a corollary of Hille's theorem. In the upcoming section, we discuss some of the important functional inequalities and their discrete versions.
\section{On Hermite--Hadamard type inequalities}
In the classical function theory, there are several crucial functional inequalities. One of the most important inequalities is Ostrowski's inequality(see \cite{Ostrowski}). There are several versions of this inequality. The classical one can be stated as follows:

"If $f\in L[a,b]$ is a differentiable mapping in $]a,b[$ such that $||f'||\leq \M$; then for any $x\in[a,b]$ the following inequality is satisfied
\Eq{*}{
\left |f(x)-\dfrac{1}{b-a}\int_{b}^{a}f(t)dt\right|\leq \M(b-a)\left[\dfrac{1}{4}+\dfrac{\left(x-\dfrac{a+b}{2}\right)^2}{(b-a)^2}\right]^"
.}
In short for the function $f$; Ostrowski's inequality provides the sharp upper bound for the difference of functional value at any point with the integral mean of the function $f$ or 
vice-versa. We can implement this idea in the case of sequences. We assume arithmetic mean, $(u_1+\cdots+u_n)/n$ of the sequence $\big<u_i\big>_{i=1}^{n}$ is denoted by $\M_{u_n}$. The height of the sequence can be defined as follows
\Eq{*}
{
\H:=\sup_{i\neq j}|u_i-u_j| \quad \mbox{where} \quad i,j
\in\{1,\cdots, n\}
.}
Then for any $u_i\in\big<u_i\big>_{i=1}^{n}$, the following inequality holds
\Eq{*}{
\left|\M_{u_n}-u_i\right|\leq \H
\quad\mbox{where} \quad i\in\{1,\cdots,n\}.}
It is easily verifiable and hence the proof is not elaborated. \\

The other important inequality is the Hermite--Hadamard inequality (see \cite { Hadamard, Hermite}). It can be formulated as:

"If $f:[a,b]\to\R$ is a convex function, then the following functional inequality is always satisfied
\Eq{*}{
f\left(\dfrac{a+b}{2}\right)\leq\int_a^{b}f(t)dt\leq \dfrac{f(a)+f(b)}{2}^".
}
In other words, the basic version of Hermite--Hadamard inequality gives bounds to the integral mean of a convex function defined in a compact interval of $\R$. Since then many versions of this inequality have been investigated in various function classes. A Hermite--Hadamard type inequality for the subadditive function was introduced in the paper \cite{Ali}. It is showed that if $f:[0,\infty[\to\R$ is a subadditive function, then for any $0<a<b<\infty$ the following inequality holds
\Eq{*}{
\dfrac{1}{2}f(a+b)\leq \int_a^{b}f(t)dt\leq \dfrac{1}{a}\int_0^{a}f(t)dt+\dfrac{1}{b}\int_0^{b}f(t)dt.
}
Motivated by this result, we attempt to formulate a result that gives upper and lower bounds to the arithmetic mean of a finite subadditive sequence.
This result can be treated as the Hermite-- Hadamard type version for subadditive sequences. 
\Thm{100}
{Let $\big<u_i\big>_{i=1}^{n}$ be a subadditive sequence. Then the following inequalities hold
\Eq{100}{ \dfrac{1}{2}\left(1+\dfrac{2}{n}\right)u_n-\dfrac{1}{n}u_{n/2}\leq &\M_{u_n}\leq\dfrac{n+1}{2}u_1\quad\quad\mbox{if $n$ is even}\\
&\mbox{and}\\
\dfrac{1}{2}\left(1+\dfrac{1}{n}\right)u_n\leq &\M_{u_n}\leq\dfrac{n+1}{2}u_1 \quad\quad\mbox{if $n$ is odd}
;}
where $\M_{u_n}$ is used to denote the arithmetic mean of the sequence.
}
\begin{proof}
Since the sequence $\big<u_i\big>_{i=1}^{n}$ is subadditive, the following system of inequalities hold
\Eq{*}{
u_1=1u_1,\quad
u_2\leq u_1+u_1= 2u_1,
&\cdots \cdots\cdots,\mbox{and}\quad
u_n\leq \sum_{i=1}^{n}u_1= n u_1
.}
Summing up all these inequalities side by side, we obtain
\Eq{*}{
u_1+\cdots+u_n\leq (1+\cdots+n)u_1=\dfrac{n(n+1)}{2}u_1.
}
By dividing both sides by $n$, we arrive at the rightmost sides of inequalities in \eq{100}. \\

Now to show the leftmost inequalities, we consider two cases. First, we assume that $n$ is even. Then we can construct the following system of $n/2$ inequalities
\Eq{*}{
u_n\leq u_1+u_{n-1},\quad
u_n\leq u_2+u_{n-2},
\cdots\cdots\cdots,\mbox{and}\quad
u_n\leq u_{_{n/2}}+u_{_{n/2}}.
}
Multiplying both sides of these inequalities by $2$ and then summing up all side by side, we get
\Eq{*}{
(n+2)u_n\leq 2(u_1+\cdots +u_n)+2u_{n/2}.
}
After dividing both sides of this inequality by $2n$ and rearranging the terms we arrive at 
\Eq{*}{
\dfrac{1}{2}\left(1+\dfrac{2}{n}\right)u_n-\dfrac{1}{n}u_{n/2}\leq \dfrac{u_1+\cdots u_n}{n}.
}
Which is the leftmost part in the first inequality of \eq{100}. \\

Again assuming $n\in\N$ is odd; we can formulate the following inequities 
\Eq{*}{
u_n\leq u_1+u_{n-1}\quad
u_n\leq u_2+u_{n-2},
\cdots \cdots\cdots,\mbox{and}\quad
u_n\leq u_{_{(n-1)/2}}+u_{_{(n+1)/2}}.
}
As before, adding up all the inequalities of the  above system, we get
\Eq{*}{
\dfrac{n-1}{2}u_n\leq u_1+\cdots+u_{n-1}.
}
Adding $u_n$ on both sides of the above inequality and then dividing by $n$; 
we obtain the leftmost of the second inequality in \eq{100}.\\

Summarizing both cases, we establish that \eq{100} holds for any $n\in\N$.
\end{proof}
\section{On Approximately periodic Sequences}
Let $\varepsilon>0$ and $\L\in\N$ such that $1\leq\L<n.$
A sequence $\big<u_i\big>_{i=1}^{\infty}$ is said to be $\varepsilon$-periodic, if it satisfies the following discrete functional inequality 
\Eq{104}{
|u_{_{i+k\L}}-u_{{_i}}|\leq\varepsilon \quad\mbox{for all}\qquad k\in\Z\qquad\mbox{with}\qquad 
i+k\L\in\N;}
where $\L$ is termed as the period of the sequence.  \\

The theorem below can be treated as a stability theorem for approximately periodic sequences. Here we will present the result for the finite version of $\varepsilon$-periodicity. But one can generalize this decomposition result for a $\varepsilon$-periodic sequence with infinite elements. 
\Thm{103}{
Assume $\big<u_i\big>_{i=1}^{n}$ be a $\varepsilon$-periodic sequence with the period $\L$. Then there exists a periodic sequence $\big<v_i\big>_{i=1}^{n}$ satisfying the inequality
\Eq{105}{
\left|u_i-v_i\right|\leq\dfrac{\varepsilon}{2}\qquad\mbox{for all}\qquad i\in\{1,\cdots,n\}
.}
Conversely, suppose \eq{105} holds, where $\big<v_i\big>_{i=1}^{n}$ is a periodic sequence with periodicity $\L$. Then $\big<u_i\big>_{i=1}^{n}$ is a $\varepsilon$-periodic sequence.
}
\begin{proof}
Since $\big<u_i\big>_{i=1}^{n}$ be a $\varepsilon$-periodic, it will satisfy \eq{104}. Since the sequence is finite, it is bounded as well. We construct the sequence $\big<v_i\big>_{i=1}^{n}$ as follows
\Eq{*}{
v_i:=\dfrac{{\max}{\Big(u_{_{i+k\L}}\Big)}+\min{\Big(u_{_{i+k\L}}\Big)}}{2}\quad \mbox{for all} \quad k\in\Z \quad \mbox{such that}\quad i+k\L\in\{1,\cdots,n\}.
}
We assume an arbitrary $k'\in\Z$  such that $i+k'\L\in\{1,\cdots,n\}.$ Then by using the definition of the sequence $\big<v_i\big>_{i=1}^{n}$, we can compute the following equality
\Eq{*}{
v_{_{i+k'\L}}=\dfrac{{\max}{\Big(u_{_{(i+k'\L)+k\L}}\Big)}+\min{\Big(u_{_{(i+k'\L)+k\L}}\Big)}}{2}
=\dfrac{{\max}{\Big(u_{_{i+k\L}}\Big)}+\min{\Big(u_{_{i+k\L}}\Big)}}{2}=v_i.
}
This shows that $\big<v_i\big>_{i=1}^{n}$ possesses periodicity with the period $\L$.\\

To establish the second part of the first assertion; we utilize \eq{104} and obtain the following two inequalities
\Eq{*}{
u_i-v_i&\leq u_i-\dfrac{{\max}{\Big(u_{_{i+k\L}}\Big)}+\min{\Big(u_{_{i+k\L}}\Big)}}{2}\leq \dfrac{u_i-\min{\Big(u_{_{i+k\L}}\Big)}}{2}\leq\dfrac{\varepsilon}{2}\\
&\qquad\qquad\quad\qquad\qquad \qquad\mbox{and}\\
v_i-u_i&\leq \dfrac{{\max}{\Big(u_{_{i+k\L}}\Big)}+\min{\Big(u_{_{i+k\L}}\Big)}}{2}-u_i\leq \dfrac{\max{\Big(u_{_{i+k\L}}\Big)}-u_i}{2}\leq\dfrac{\varepsilon}{2}.
}
This two inequalities together yields $||u_i-v_i||\leq\dfrac{\varepsilon}{2}$ and establishes the first assertion.\\

To show the reverse part, we consider an arbitrary $i\in\{1,\cdots,n\}$. Let $k'\in\Z$ such that $i+k'\L\in\{1,\cdots,n\}$, we have that both $u_i$, $u_{_{i+k'\L}}\in \big<u_i\big>_{i=1}^{n}$. Then by using \eq{105} and periodicity of the sequence $\big<v_i\big>_{i=1}^{n}$, we can proceed as follows
\Eq{*}{
\left|u_i-u_{_{i+k'\L}}\right|
&=\left|u_i-v_i+v_{{i+k'\L}}-u_{_{i+k'\L}}\right|\\
&\leq \left|u_i-v_i\right|+\left|u_{{i+k'\L}}-v_{_{i+k'\L}}\right|\\
&\leq \dfrac{\varepsilon}{2}+\dfrac{\varepsilon}{2}=\varepsilon.
}
Since $i\in\N$ and $k'\in\Z$ are arbitrary, this yields \eq{104}. In other words the sequence $\big<u_i\big>_{i=1}^{n}$ possesses $\varepsilon$-periodicity with the periof $\L$. This completes the proof.
\end{proof}
In \eq{104} if $\varepsilon=0$, then $\big<u_i\big>_{i=1}^{n}$ is nothing but a ordinary periodic sequence with the period $\L$. A constant sequence is periodic with periodicity, $\L=1$. Conversely,  a periodic sequence with unit periodicity implies a constant sequence. For this reason, in the study of ordinary periodic sequences, we are only interested in those containing periodicity of at least $2$.
\Prp{101}{Let $\big<u_n\big>_{n=1}^{\infty}$ be a periodic sequence with period $\L$. Then there exists a continuous periodic function  
$f:[1,\infty[\to\R$ that possesses the same period and satisfies $f(n)=u_n$ for all $n\in\N.$}
\begin{proof}
To prove this theorem, we construct the function $f:[1,\infty[\to\R$ as mentioned in \eq{4000}.
\newline
The function $f$ is continuous.
By replacing $u=n$ (or by taking $t=1$); it is evident that $f(n)=u_n$ holds for all $\N$. Thus we are only left to show that $f$ possesses periodicity with the periodic length $\L$.
We consider an arbitrary $u\in [n,n+1] \subset[1,\infty[$ $\Big(n\in\N\Big)$ and $k\in\Z$ such that $u+k\L\in [1,\infty[.$ It yields $u+k\L\in[n+k\L,n+1+k\L]$ and allows us to compute the following extended equality
\Eq{*}{
f(u+k\L)&=t\big(u_{_{n+k\L}}\big)+(1-t)\big(u_{_{n+1+k\L}}\big)
=tu_n+(1-t)u_{n+1}
=f(u).
}
It establishes that $f:[1,\infty[\to\R$ is a periodic function having the periodic length $\L$. This completes the proof.
\end{proof}
\section{On Characterizations of Periodic Sequences}
Let $\big<u_n\big>_{n=1}^{\infty}$ be a sequence. For $n\in\N\cup\{0\}$, we define the terms $S_0=0$ and $S_n:=u_1+\cdots+u_n$. Then the sequence $\big<S_n\big>_{n=0}^{\infty}$ is said to be the sequence of partial sums due to the sequence $\big<u_n\big>_{n=1}^{\infty}$. In the very first result of this section, we will see in the case of a periodic sequence, how its sequence of partial sums correlates.
 
\Thm{500}{If $\big<u_n\big>_{n=1}^{\infty}$ be a periodic sequence with period $\L$. 
Then for any  $i\in\{1,\cdots,\L\}$, the sequence $\big<v_{_n}^{^{i}}\big>_{n=1}^{\infty}$ with  $v_{_{n}}^{^{i}}=S_{_{i+(n-1)\L}}-S_{_{(n-1)\L}}$  is a constant sequence.
\\
On the other hand, If for each of the $i\in\{1,\cdots,\L\}$, all of the $L$ sequences 
$\big<v_{_n}^{^{i}}\big>_{n=1}^{\infty}$ are constants, then this will imply $\big<u_n\big>_{n=1}^{\infty}$ is a periodic sequence with periodicity $\L$.
}
\begin{proof}
To prove the first assertion, we choose $i\in\{1,\cdots,\L\}$ arbitrarily and compute each element of the sequence $\big<v_{_n}^{^{i}}\big>_{n=1}^{\infty}$ as follows 
\Eq{*}{
\mbox {For} \qquad n=1 \qquad v_{_1}^{^{i}}=S_i-S_0&=u_1+\cdots+u_i.\\
\mbox {And for all other $n\in\N$ $(n>1)$} \qquad v_{_n}^{^{i}}&=S_{_{i+(n-1)\L}}-S_{_{(n-1){\L}}}\\
&=u_{_{1+(n-1)\L}}+\cdots+u_{_{i+(n-1)\L}}\\
&=u_1+\cdots+u_i.
}
This shows that $\big<v_{_n}^{^{i}}\big>_{n=1}^{\infty}$ is a constant sequence. \\

To establish the reverse implication, we choose  $i\in\{1,\cdots,\L\}$ arbitrarily. Now we can compute the following two equalities
\Eq{10}{
v_{_1}^{^{i}}&=S_i-S_0=u_1+\cdots+u_i
\\
 \qquad\qquad\qquad\mbox{and for all other $n\in\N$ $(n>1)$, we have}\\ 
v_{_{n}}^{^{i}}&=S_{_{i+(n-1)\L}}-S_{_{(n-1)\L}}=u_{_{1+(n-1)\L}}+\cdots+u_{_{i+(n-1)\L}}.
}
Now we consider two cases. At first, by replacing $i=1$ in \eq{10}, we have
\Eq{*}{
v_{_1}^{^{1}}=u_1\qquad\mbox{and}\qquad v_{_{n}}^{^{1}}=u_{_{1+(n-1)\L}} \quad \mbox{for all other $n\in\N$ $(n>1)$.}
}
Since $\big<v_{_n}^{^{i}}\big>_{n=0}^{\infty}$ is  constant, the two equations above together yields $u_1=u_{_{1+n\L}}$ for all $n\in\N$.\\

In the second case, considering any $i\in\{2,\cdots,\L\}$ and by utilizing \eq{10}, we can compute the following
\Eq{290}{
v_{_1}^{^{i}}-v_{_1}^{^{i-1}}=u_i\qquad\mbox{and}\qquad v_{_n}^{^{i}}-v_{_n}^{^{i-1}}=u_{_{i+(n-1)\L}} \quad \mbox{for all other $n\in\N$ $(n>1)$}.
}
Since the difference of two constant sequences $\big<v_{_n}^{^{i}}\big>_{n=1}^{\infty}$ and 
$\big<v_{_n}^{^{i-1}}\big>_{n=1}^{\infty}$    results in the constant sequence $\Big<v_{_n}^{^{i}}-v_{_n}^{^{i-1}}\Big>_{n=1}^{\infty}$;  the equations in \eq{290} implies $u_i=u_{i+n\L}$ for all $n\in\N$.\\

Together these two cases establish that $\big<u_n\big>_{n=1}^{\infty}$ is a constant sequence and completes the proof.
\end{proof}
In the second result, we will prove that any periodic sequence can be expressed as a combination of constant sequences. But before that, we need to go through three concepts. Namely distinctness, merger, and inclusion of sequences. One can see that these definitions are motivated by basic set-theoretic concepts.\\

Two sequences $\big<u_n\big>_{n=1}^{\infty}$ and $\big<v_n\big>_{n=1}^{\infty}$ are said to be distinct if they have no elements in common. Symbolically, we will represent it as $\big<u_n\big>_{n=1}^{\infty}\bigwedge \big<v_n\big>_{n=1}^{\infty}=\phi$. \\

Let $\L\in\N$ be fixed. Then by Euclid's division lemma for any $n\in\N$, we have the expression $n=i+k\L$; where $k\in\N\cup\{0\}$ and $i\in{1,\cdots,\L}$ is unique. Now for any sequence $\big<u_n\big>_{n=1}^{\infty}$; we can extract $\L$ number of sub-sequences as $\Big<v^{^{i}}_{_{k}}\Big>_{k=0}^{\infty}$ ($i\in\{1,\cdots,\L\}$). Clearly for any  $u_n\in\big<u_n\big>_{n=1}^{\infty}$, we have $u_n=v^{^{i}}_{_{k}}$ and this representation is exclusive. This also shows that any element that is chosen from one of the sub-sequences$\Big<v^{^{i}}_{_{k}}\Big>_{k=0}^{\infty}$'s; we can determine the exact position of it in the original sequence $\big<u_n\big>_{n=1}^{\infty}$ as well. In other words, element-wise ordered merger of these $\L$ sub-sequences $\Big<v^{^{i}}_{_{k}}\Big>_{k=0}^{\infty}$'s will eventually bring back the sequence $\big<u_n\big>_{n=1}^{\infty}$. Mathematically, we are going to denote it as  $\overset{i=\L}{\underset{i=1}{\bigvee}}\Big<v^{^{i}}_{_{k}}\Big>_{k=0}^{\infty}=\big<u_n\big>_{n=1}^{\infty}.$\\

If all the elements of the sequence $\big<v_n\big>_{n=1}^{\infty}$ are also elements of the sequence $\big<u_n\big>_{n=1}^{\infty}$; we will say $\big<v_n\big>_{n=1}^{\infty}$ is contained in the sequence $\big<u_n\big>_{n=1}^{\infty}$. We represent this inclusion property by $\big<v_n\big>_{n=1}^{\infty}\sqsubseteq\big<u_n\big>_{n=1}^{\infty}.$
\Thm{102}{
Let $\big<u_n\big>_{n=1}^{\infty}$ be a periodic sequence with the period $\L$. Then there can exists $\ell$ $(\ell\leq\L)$ number of constant sub-sequences $\Big<v^{^{i}}_{_{k}}\Big>_{k=0}^{\infty}$ 
$(i=1,\cdots,\ell)$  such that the following two conditions are satisfied
\Eq{102}{(i)\quad \mbox{If}\quad  i\neq j, \quad \mbox{then} \quad \Big<v^{^{i}}_{_{k}}\Big>_{k=0}^{\infty} \bigwedge\Big<v^{^{j}}_{_{k}}\Big>_{k=0}^{\infty}=\phi
\quad 
\mbox{and}
\quad 
(ii)\quad \overset{i=\ell}{\underset{i=1}{\bigvee}}\Big<v^{^{i}}_{_{k}}\Big>_{k=0}^{\infty}=\big<u_n\big>_{n=1}^{\infty}. 
}
Conversely, if $\big<u_n\big>_{n=1}^{\infty}$ be a sequence having a total $\L$ number of constant sub-sequences $\Big<v^{^{i}}_{_{k}}\Big>_{k=0}^{\infty}$  $(i\in\{1,\cdots,\L\})$ satisfying the two conditions mentioned in \eq{102}. Then $\big<u_n\big>_{n=1}^{\infty}$ is a periodic sequence with period $\L$.
}
\begin{proof}
If $\big<u_n\big>_{n=1}^{\infty}$ itself is a constant sequence there is nothing to show.
To show this theorem, we assume that $\L>1$ and we construct total $\L$ number of subsequences $\Big<v^{^{i}}_{_{k}}\Big>_{k=0}^{\infty}$ $\Big( i\in\{1,\cdots,\L\}\Big)$ as follows
\Eq{*}{
\Big<v^{^{1}}_{_{k}}\Big>_{k=0}^{\infty}:=\Big<u_{_{1+k\L}}\Big>_{k=0}^{\infty},
\cdots \cdots,\mbox{and}\quad
\Big<v^{^{L}}_{_{k}}\Big>_{k=0}^{\infty}:=\Big<u_{_{(k+1)\L}}\Big>_{k=0}^{\infty}\qquad (k\in\N).
}
Since, $\big<u_n\big>_{n=1}^{\infty}$ is a periodic with the period $\L$; this yields that 
$\Big<v^{^{i}}_{_{k}}\Big>_{k=0}^{\infty}$ $\Big( i\in\{1,\cdots,\L\}\Big)$ are constant sub-sequences. Without loss of generality, we assume that all these sequences are unique. Suppose If not, we consider the merger of two (or more) similar sequences $\Big<v^{^{i}}_{_{k}}\Big>_{k=0}^{\infty}$ and $\Big<v^{^{j}}_{_{k}}\Big>_{k=0}^{\infty}$ as a single sequence $\Big<v^{^{i\vee j}}_{_{k}}\Big>_{k=0}^{\infty}$; mathematically  
$\Big<v^{^{i\vee j}}_{_{k}}\Big>_{k=0}^{\infty}:=\Big<v^{^{i}}_{_{k}}\Big>_{k=0}^{\infty}\bigvee \Big<v^{^{j}}_{_{k}}\Big>_{k=0}^{\infty}$. This will ensure that the total number of constant sub-sequences remains less than $\L$. And the establishment of this proposition can be carried out similarly. \\

The first assertion is straightforward. Since each of the sub-sequences $\Big<v^{^{i}}_{_{k}}\Big>_{k=0}^{\infty}$ $i\in\{1,\cdots\ell\}$ are constants; all of these must be distinct. This is the first equality of \eq{102}\\

Now we will show the second part of \eq{102}.
$\quad \overset{i=\L}{\underset{i=1}{\bigvee}}\Big<v^{^{i}}_{_{k}}\Big>_{k=0}^{\infty}\sqsubseteq\big<u_n\big>_{n=1}^{\infty}$ is obvious.
To show the reverse inclusion, We assume $u_n\in \big<u_n\big>_{n=1}^{\infty}$. For $n$ and $\L\in\N$, there exist a specific $k\in\N\cup\{0\}$ and a unique $i\in\{1,\cdots,\L\}$ such that $n=i+k\L$ holds. This indicates the following
\Eq{*}{
u_n=u_{_{i+k\L}}\in \Big<v^{^{i}}_{_{k}}\Big>_{k=0}^{\infty}\in \quad \overset{i=\L}{\underset{i=1}{\bigvee}}\Big<v^{^{i}}_{_{k}}\Big>_{k=0}^{\infty} \quad(i\in \{1,\cdots \L\})
.}
Hence $\big<u_n\big>_{n=1}^{\infty}\sqsubseteq \overset{i=\L}{\underset{i=1}{\bigvee}}\Big<v^{^{i}}_{_{k}}\Big>_{k=0}^{\infty}$ holds and this establishes the second implication of \eq{102}.\\

To show the reverse part, we assume that $u_n\in\big<u_n\big>_{n=1}^{\infty}$ be arbitrary.  The two conditions mentioned in \eq{102} together yields that out of total $\L$ number of distinct constant sub-sequences there must exists a sub-sequence $\Big<v^{^{i}}_{_{k}}\Big>_{k=0}^{\infty}$ ($i\in\{1,\cdots,\L\}$) satisfying the following equality
\Eq{98}{
v^{^{i}}_{_{k}}=u_{_{i+k\L}}=u_n \quad \quad \quad (n\in\N).
}
Without loss of generality, we choose $k'\in\Z$  such that $n+k'\L\in\N$. Since $\Big<v^{^{i}}_{_{k}}\Big>_{k=0}^{\infty}$ is constant; by utilizing \eq{98}, we can get the following extended equality
\Eq{*}{
u_{_{n+k'\L}}=u_{_{(i+k\L)+k'\L}}=u_{_{i+(k+k')\L}}=v^{^{i}}_{_{k+k'}}=v^{^{i}}_{_{k}}=u_n \quad \mbox{for all} \quad n\in\N.
}
Since $k'$ is arbitrary, the above equality establishes that the sequence $\big<u_n\big>_{n=1}^{\infty}$ possesses periodicity of length $\L$. And it completes the proof of the statement.
\end{proof}
One can investigate the concepts of 
$\varepsilon$-periodicity and subadditivity in stochastic processes, potentially leading to new stability theorems. Results such as the interpolation of subadditive sequences via continuous subadditive functions and the decomposition of 
$\varepsilon$-periodicity could be valuable in addressing various mathematical modeling 
problems. These approaches may provide novel solutions in fields like optimization, numerical analysis, dynamical systems, and computational mathematics.\\

\textbf{Conflicts of Interest Statement:}
The author declare that there are no conflicts of interest related to the research, authorship, and/or publication of this manuscript. The research presented is an original work and has not been submitted elsewhere for publication.\\

\textbf{Data Availability Statement:}
This manuscript does not utilize any datasets. All mathematical derivations, theorems, and results presented are derived analytically. No external data was generated or analyzed during the study.


\begin{thebibliography}{1}
\bibitem{Fekete}
Fekete, Michael. 
\newblock {"Über die Verteilung der Wurzeln bei gewissen algebraischen Gleichungen mit ganzzahligen Koeffizienten."}
\newblock { Mathematische Zeitschrift}
17, no. 1 (1923): 228-249.
\bibitem{Hillee}
Hille, Einar. 
\newblock {"Non-oscillation theorems."}
\newblock {Transactions of the American Mathematical Society}
 64, no. 2 (1948): 234-252. 
 
 \bibitem{Hille}
Hille, Einar, and Ralph S. Phillips. 
\newblock {Functional analysis and semi-groups.}
\newblock { Vol. 31. New York: American Mathematical Society,}
 1948.  

\bibitem{Schechter}
Schechter, Eric. 
\newblock {Handbook of Analysis and its Foundations.}
\newblock {Academic Press,}
 1996.
 
\bibitem{Matkowski}
Matkowski, Janusz, and Tadeusz Świątkowski. 
\newblock {"On subadditive functions."}
\newblock {Proceedings of the American Mathematical Society}
119, no. 1 (1993): 187-197. 

\bibitem{Rosenbaum}
Rosenbaum, Robert A. 
\newblock {"Sub-additive functions."}
(1950): 227-247.

\bibitem{Matkowskii}
Matkowski, Janusz. 
\newblock {"Subadditive periodic functions."}
\newblock { Opuscula Mathematica}
 31, no. 1 (2011): 75-96.
 
\bibitem{Goswamii} 
Goswami, Angshuman R. 
\newblock {"Generalization of Subadditive, Monotone and Convex Functions."}
\newblock {arXiv preprint}
arXiv:2308.00704 (2023).

\bibitem{Burai}
Burai, Pál, and Árpád Száz. 
\newblock {"Homogeneity properties of subadditive functions."}
\newblock {In Annales Mathematicae et Informaticae,}
 vol. 32, pp. 189-201. 
Eszterházy Károly College, Institute of Mathematics and Computer Science, 2005.

\bibitem{Arpad}
Burai, Pál, and Arp\'ad Száz. 
\newblock {"Relationships between homogeneity, subadditivity and convexity properties."} \newblock {Publikacije Elektrotehničkog fakulteta.} 
Serija Matematika
 2005: 77-87.

\bibitem{ Kominek} 
 Kominek, Zygfryd. 
\newblock {"On approximately subadditive functions."}
\newblock {Demonstratio Mathematica}
23, no. 1 (1990): 155-160.
 
\bibitem{Ostrowski}
Ostrowski, Alexander. 
\newblock {"Über die Absolutabweichung einer differentiierbaren Funktion von ihrem Integralmittelwert."}
\newblock{Commentarii Mathematici Helvetici} 
10, no. 1 (1937): 226-227.

\bibitem{Hadamard}
J. Hadamard.
\newblock{Étude sur les propriétés des fonctions entières et en particulier d'une fonction considérée par Riemann.}
\newblock{\em Journal de mathématiques pures et appliquées}
 9 (1893): 171-215.

\bibitem{Hermite}
Ch. Hermite.
\newblock{Sur deux limites d’une intégrale définie}.
\newblock{\em Mathesis}
3(1) (1883): 1-82.

 \bibitem{Ali}
Sarikaya, Mehmet Zeki, and Muhammad Aamir Ali. \newblock {"Hermite-Hadamard type inequalities and related inequalities for subadditive functions."}
\newblock { Miskolc Mathematical Notes}
22, no. 2 (2021): 929-937.

\bibitem{Hyers}
Hyers, Donald H. 
\newblock{"On the stability of the linear functional equation."}
\newblock{Proceedings of the National Academy of Sciences}
27, no. 4 (1941): 222-224.

\bibitem{Ulam}
Hyers, Donald H., George Isac, and Themistocles M. Rassias. 
\newblock{"Approximately convex functions."}
\newblock{In Stability of Functional Equations in Several Variables,}
pp. 166-179. Boston, MA: Birkhäuser Boston, 1952.

\bibitem{Gilbert}
Gilbert, Edgard N., and John Riordan. 
\newblock{"Symmetry types of periodic sequences."} 
\newblock{Illinois Journal of Mathematics}
 5, no. 4 (1961): 657-665.
 
\bibitem{Heuberger} 
Heuberger, Clemens, and Daniel Krenn. 
\newblock{"Asymptotic analysis of regular sequences."}
\newblock{Algorithmica}
 82, no. 3 (2020): 429-508.
 
\bibitem{Muchnik}
Muchnik, An, A. Semenov, and M. Ushakov. 
\newblock{"Almost periodic sequences."}
\newblock{ Theoretical Computer Science}
304 (2003): 1-33.

\bibitem{Kurshan}
Kurshan, R. P., and B. Gopinath. 
\newblock{"Recursively generated periodic sequences."}
\newblock{Canadian Journal of Mathematics}
26, no. 6 (1974): 1356-1371.

\bibitem{Khamidullin}
Kel'Manov, A. V., and S. A. Khamidullin. \newblock{"Subsequences in a Quasi-periodic Sequence."}
\newblock{Computational, Mathematics and Mathematical Physics}
 4, no. 5 (2001).
 
\bibitem{Goswami}
Goswami, Angshuman Robin. 
\newblock{"Decomposition of Approximately Monotone and Convex Sequences."}
arXiv preprint arXiv:2404.15362 (2024).
\end{thebibliography}
\end{document}